\ifx\LocalTwelvePtArticleFormatLoaded\undefined
  \documentclass[a4paper,12pt]{article}     
  \usepackage{array}                        
  \usepackage{theorem}                      
  \usepackage{amsmath,amscd,amssymb}        
  \usepackage{latexsym}                     
  \usepackage[german,english]{babel}        
  \usepackage[applemac]{inputenc}           
\fi

\newcommand{\SL}{\mathrm {SL}}
\newcommand{\PSL}{\mathrm {PSL}}

\newcommand{\SP}{\mathrm {Sp}}

\newcommand{\Coker}{\mathrm{Coker}}
\newcommand{\Ker}{\mathrm{Ker}}
\newcommand{\Gr}{\mathrm{Gr}}

\newcommand{\into}{\hookrightarrow}

\newcommand{\et}{\mbox{\scriptsize \'et}}
\newcommand{\Het}{H_{\et}}

\newcommand{\cA}{\mathcal {A}}
\newcommand{\cL}{\mathcal {L}}
\newcommand{\Gal}{\mathrm{Gal}}
\newcommand{\cZ}{\mathcal {Z}}

\newcommand{\cI}{\mathcal {I}}

\newcommand{\hX}{\widehat{X}}

\newcommand{\cP}{\mathcal{P}}

\newcommand{\cO}{\mathcal{O}}

\newcommand{\ZZ}{\mathbb Z}
\newcommand{\QQ}{\mathbb Q}
\newcommand{\bQQ}{\overline{\QQ}}

\newcommand{\AAA}{\mathbb A}
\newcommand{\CC}{\mathbb C}
\newcommand{\FF}{\mathbb F}
\newcommand{\PP}{\mathbb P}

\newcommand{\Aut}{\mathrm{Aut}}
\newcommand{\Frob}{\mathrm{Frob}}
\newcommand{\Tr}{\mathrm{tr}}
\newcommand{\tY}{\widetilde{Y}}
\newcommand{\tU}{\widetilde{U}}

\newcommand{\tW}{\widehat{W}}

\newcommand{\tN}{\widetilde{N}}
\newcommand{\tNn}{{\tN}^{(0)}}
\newcommand{\tNe}{{\tN}^{(1)}}
\newcommand{\tNt}{{\tN}^{(2)}}
\newcommand{\tNd}{{\tN}^{(3)}}
\newcommand{\Nn}{N^{(0)}}
\newcommand{\Ne}{N^{(1)}}
\newcommand{\Nt}{N^{(2)}}
\newcommand{\Nd}{N^{(3)}}

\newcommand{\hPP}{\widehat{\PP}}

\newcommand{\tPP}{\widetilde{\PP}}
\newcommand{\PPn}{\PP^{(0)}}
\newcommand{\PPe}{\PP^{(1)}}
\newcommand{\PPt}{\PP^{(2)}}
\newcommand{\PPd}{\PP^{(3)}}
\newcommand{\tPPn}{\tPP^{(0)}}
\newcommand{\tPPe}{\tPP^{(1)}}
\newcommand{\tPPt}{\tPP^{(2)}}
\newcommand{\tPPd}{\tPP^{(3)}}

\newcommand{\tp}{\tilde{p}}
\newcommand{\tq}{\tilde{q}}
\newcommand{\hq}{\hat{q}}

\newcommand{\tf}{\tilde{f}}

\newcommand{\Ce}{C^{(1)}}
\newcommand{\Ct}{C^{(2)}}
\newcommand{\Cd}{C^{(3)}}
\newcommand{\tCe}{\widetilde{C}^{(1)}}
\newcommand{\tCt}{\widetilde{C}^{(2)}}
\newcommand{\tCd}{\widetilde{C}^{(3)}}

\newcommand{\Sing}{\mathrm{Sing}}

\newtheorem{thm}{Theorem}[section]
\newtheorem{cor}[thm]{Corollary}
\newtheorem{lem}[thm]{Lemma}
\newtheorem{rmk}[thm]{Remark}

\newtheorem{pro}[thm]{Proposition}

\newsavebox{\pfaux}
\newenvironment{pf}[2]{\noindent {\scshape {#1}}\sbox{\pfaux}{#2}}
{\hfill \usebox{\pfaux}\medskip}

\begin{document}


\bigskip

\bigskip

\centerline{\Large\bfseries The modularity of the }
\centerline{\Large\bfseries Barth-Nieto quintic and its relatives}

\vskip 50pt

\centerline{\large K.~Hulek, J.~Spandaw, B.~van Geemen, D. van Straten}

\bigskip
\bigskip

\noindent {\bfseries Abstract.}  {\footnotesize
The moduli space of $(1,3)$-polarized
abelian surfaces with full level-2 structure
is birational to a double cover of the Barth-Nieto quintic.
Barth and Nieto have shown that these varieties have Calabi-Yau models
$Z$ and~$Y$, respectively. In this paper we apply the Weil
conjectures to show that $Y$ and~$Z$
are rigid and we prove that the $L$-function of their
common third \'etale cohomology group is modular,
as predicted by a conjecture of Fontaine and Mazur.
The corresponding modular form is the unique normalized cusp form of
weight~4 for the group~$\Gamma_1(6)$.  By Tate's conjecture,
this should imply that $Y$, the fibred square
of the universal elliptic curve~$S_1(6)$, and Verrill's
rigid Calabi-Yau~$\cZ_{A_3}$, which all have the same $L$-function,
 are in correspondence over~$\QQ$.
We show that this is indeed the case by giving explicit maps.
}

\bigskip

\section{Introduction}

The Barth-Nieto quintic is the variety given by the equations
$$
N=\left\{\sum_{i=0}^5 x_i=\sum_{i=0}^5 \frac1{x_i}=0\right\}\subset\PP_5.
$$
This singular quintic threefold in $\PP_4=\{\sum_{i=0}^5x_i=0\}$
 was studied by Barth and Nieto in \cite{BN}. They show that it
parametrizes Heisenberg $H_{2,2}$-invariant quartics in~$\PP_3$ which
contain a line.
A smooth such quartic then contains 32~lines which form two $H_{2,2}$-orbits
of 16~disjoint lines each. Taking the double cover branched
along these orbits of lines gives two abelian surfaces which are dual and
which have a polarization of type (1,3).
This defines a map $\cA_{1,3}(2)\dashrightarrow N$ of degree~2, where
$\cA_{1,3}(2)$ is the moduli space of abelian surfaces with
(1,3)-polarization and a full level-2 structure.
The  moduli space $\cA_{1,3}(2)$ is birationally equivalent to the inverse
image
$\tN$ of~$N$ under the the double cover of~$\PP_5$ branched along the union
of the 6~hyperplanes~$\{x_k=0\}$.
(All tildes \lq\lq\ $\widetilde{}$ \rq\rq\ in this paper denote double covers;
resolutions of singularities are denoted by a hat~\lq\lq\ $\widehat{}$ \rq\rq.)

The varieties $N$ and $\tN$ have smooth (strict) Calabi-Yau models,
denoted by $Y$ and $Z$ respectively. Thus they have trivial canonical
bundle  and
$h^{i}(\cO)=0$ for $0<i<3$. The Euler numbers are $e(Y)=100$ and $e(Z)=80$.
The paper \cite{GH2} gives a different proof
for the existence of the Calabi-Yau model of $\tN$,
using a Siegel modular form and the birationality of $\cA_{1,3}(2)$
and~$\tN$.

In this paper
we prove that $Y$ and~$Z$ are rigid. In fact, this determines all Hodge
numbers and we obtain the following result.

\bigskip

\noindent {\bfseries Theorem \ref{th21}}
\emph{
 Both $Y$ and~$Z$ have Hodge numbers $h^{p,q}=0$ except the following:
$$
h^{0,0}=h^{3,0}=h^{0,3}=h^{3,3}=1,
$$
$$
h^{1,1}(Y)=h^{2,2}(Y)=50,
\qquad h^{1,1}(Z)=h^{2,2}(Z)=40.
$$
In particular, both manifolds are rigid.
}

\bigskip

The most difficult part of the proof of Theorem~\ref{th21} is the
determination of $h^{2,2}$ and $h^{2,1}$.
We exploit the fact that $Y$~is defined
over~$\ZZ$ and that we can replace~$Z$
by another model, $\tY$, which is defined
over~$\ZZ$ and which has
$H^3(\tY)\cong H^3(Z)$.
Following the method pioneered in \cite{vGN}, we reduce modulo some primes
where $Y$ and~$\tY$ have good reduction and explicitly count the number of
points on the reductions with the help of a computer.
Using the Lefschetz fixed point formula and the Riemann hypothesis for
varieties over finite fields  leads to the desired result.

We summarize the maps between the various varieties
by the commutative diagram
$$
\begin{CD}
Z @<<< \tY @>>> \tN @. \dashleftarrow \cA_{1,3}(2)\\
@. @VVV @VVV\\
@. Y @>>> N
\end{CD}
$$
where the horizontal and vertical maps are generically 1:1
and~2:1, respectively.
As mentioned above, we usually work with~$\tY$ instead of~$Z$,
even though it is not Calabi-Yau, since it is defined over~$\ZZ$.

\

According to Theorem~\ref{th21} the \'etale cohomology groups
$\Het^3(Y)$ and $\Het^3(\tY)$
are 2-dimensional representations of~$\Gal(\bQQ/\QQ)$.
By a conjecture of Fontaine and Mazur \cite{FM} these representations
should be modular, which, as we shall show, is indeed the case.
To make a precise statement let $\eta(q)=q^{1/12}\prod_{n=1}^\infty (1-q^n)$
be the Dedekind $\eta$-function and define
$$
f(q):=[\eta(q)\eta(q^2)\eta(q^3)\eta(q^6)]^2.
$$
This is a (normalized) newform of weight~4
with respect to the group~$\Gamma_0(6)$.
In fact, one has $S_4(\Gamma_0(6))=\CC\cdot f$,
where $S_k(\Gamma)$ denotes the vector space
of cusp forms of weight~$k$ with respect to a
congruence group~$\Gamma$.

\bigskip

\noindent {\bfseries Theorem \ref{th32}}
\emph{
There is an equality
$$
L(\Het^3(Y),s)\circeq L(\Het^3(\tY),s)\circeq L(f,s)
$$
where $L_1\circeq L_2$ means that the Euler factors of $L_1$ and~$L_2$
coincide with the possible exception of the bad primes which are $2$ and~$3$
in this case.
}

\bigskip

\begin{rmk}\upshape The rigidity of~$Y$ and the modularity of
$L(\Het^3(\tY),s)$
were already noted in \cite[p.~864]{vS}.
The $L$-function of~$Y$ was also known to R.~Livn\'{e} (unpublished).
\end{rmk}

Using the 2:1 cover $\tY\to Y$
the proof of this theorem can be reduced to checking that
$L(\Het^3(Y),s)\circeq L(f,s)$. By a result due to Serre
it suffices to check equality of the Fourier coefficients~$a_p$ for the primes
$p=5,7,11,13,17,19,23$
and~73. (These primes are good primes which represent the elements
of~$(\ZZ/24\ZZ)^*$.)

\

Given any (normalized) newform $g$ of weight 4 with respect to $\Gamma_0(N)$
with Fourier coeffcients in $\ZZ$, there are 2-dimensional
$\ell$-adic $\Gal(\bQQ/\QQ)$ representation $\rho_{g,\ell}$ on $\QQ_\ell^2$
with the property that $L(\rho_{g,\ell})=L(g,s)$.
The existence of these Galois representations was established by Deligne,
and they occur naturally in the {\'e}tale cohomology of the fibre square
$W_N:=S_1(N)\times_{X_1(N)}S_1(N)$
 of the universal elliptic curve~$S_1(N)$ over
the modular curve~$X_1(N)$.
In our case we have the Shimura isomorphism
$S_4(\Gamma_0(6))\cong S_4(\Gamma_1(6))\cong  H^{3,0}(\tW)$
where $\tW$~is a resolution of
$$
W:=W_6=S_1(6)\times_{X_1(6)}S_1(6)
$$
(cf.~\cite{SY}).
The Galois representations $\rho_{f,\ell}$ are subrepresentations of
$\Het^3(\tW,\QQ_\ell)$.

The Tate conjecture implies that if isomorphic Galois representations
$\rho_1$, $\rho_2$ occur in
the \'etale cohomology of varieties $X_1$, $X_2$ defined over $\QQ$,
then there is a correspondence, defined over $\QQ$,
between $X_1$ and $X_2$ which induces an isomorphism between $\rho_1$ and
$\rho_2$.

We say that a smooth projective variety $X$, defined over $\QQ$,
is a relative of the Barth-Nieto quintic $N$ if the Galois representation
$\rho_{f,\ell}$
occurs in $\Het^3(X_{\bar{\QQ}})$.
Tate's conjecture implies that relatives should be in correspondence.
Since this conjecture is still very much open, we tried (and succeeded!)
in establishing some correspondences.

The relative $\tY$ is birationally equivalent to $\tN$, the double cover of
$N$,
and thus there are correspondences between these relatives.
Using the explicit equations for the varieties involved
we also find a correspondence between $\tY$ and~$\tW$:

\bigskip

\noindent {\bfseries Theorem \ref{th41}}
\emph{
There exists a birational equivalence between
$S_1(6)\times_{X_1(6)}S_1(6)$ and~$Y$ which is defined over~$\QQ$.
}

\bigskip

In her investigations of rigid Calabi-Yau threefolds
H.~Verrill met another relative of the Barth-Nieto quintic.
It is a threefold denoted by $\cZ_{A_3}$ and it is constructed using the
root system $A_3$.
We establish a correspondence with this relative:

\bigskip

\noindent {\bfseries Theorem (cf.~Theorem~\ref{th43})}
\emph{
There is a birational equivalence between $Y$ and~$\cZ_{A_3}$
which is defined over~$\QQ$.
}

\bigskip

Combining this theorem with Theorems \ref{th41} gives a birational equivalence
between $\cZ_{A_3}$ and $S_1(6)\times_{X_1(6)}S_1(6)$.
Such an equivalence was first found by M.~Saito and N.~Yui
(see \cite{SY}),
but their equivalence is different from ours.

\bigskip

\noindent {\em Acknowledgement.} We thank R.~Livn\'{e} and Ch.~Schoen
for helpful comments.

\section{The two protagonists}\label{sec1}

In this section we recall the construction of
our two protagonists $Y$ and~$Z$
from \cite{BN}.

The construction of the Calabi Yau varieties
$Y$ and~$Z$ from the Nieto quintic~$N$
is summarized by the following commutative
diagram:
$$
\begin{CD}
 Z @<{\hq}<< \tY @>{\tq}>> \tNd @>{\tq_3}>> \tNt @>{\tq_2}>> \tNe
@>{\tq_1}>> \tNn\\
@. @VV{\nu}V @VV{\nu_3}V @VV{\nu_2}V @VV{\nu_1}V @VV{\nu_0}V\\
@. Y @>>{q}> \Nd @>{\sim}>{q_3}> \Nt @>>{q_2}> \Ne @>>{q_1}> \Nn:=N.
\end{CD}
$$
The three squares on the right are derived from the analogous diagram
$$
\hskip 100pt\begin{CD}
\tPPd @>{\tp_3}>> \tPPt @>{\tp_2}>> \tPPe @>{\tp_1}>> \tPPn\\
@VV{\pi_3}V @VV{\pi_2}V @VV{\pi_1}V @VV{\pi_0}V\\
\PPd @>>{p_3}> \PPt @>>{p_2}> \PPe @>>{p_1}> \PPn\cong\PP_4.
\end{CD}
$$
The varieties appearing in the first (second) diagram
have dimension three (four).
The horizontal maps are birational,
the vertical maps are generically~2:1.

We now describe the diagrams in detail.
We work in
$$
\PPn:=\{x\in\PP_5\colon x_0+\cdots+x_5=0\}\cong\PP_4.
$$
The Nieto quintic is the irreducible hypersurface
$$
\Nn:=N:=\{x\in\PPn\colon \sigma_5(x)=0\},
$$
where
$$
\sigma_5(x):=\sum_{j=0}^5 x_0\cdots \hat{x}_i\cdots x_5
=x_0\cdots x_5\sum_{j=0}^5\frac1{x_j}.
$$
Its singular locus consists of the 20~lines
$$
L_{k\ell m}:=\{x_k=x_\ell=x_m=0\}\qquad (0\le k<\ell<m\le 5)
$$
in~$\PPn$
and the 10~\lq\lq Segre\rq\rq\ points
$$
\mbox{$(1:1:1:-1:-1:-1)$ $+$ permutations.}
$$
These Segre points are ordinary double points of~$N$.
The lines~$L_{k\ell m}$ intersect in the 15~points
$$
P_{k\ell m n}=\{x_k=x_\ell=x_m=x_n=0\} \qquad (0\le k<\ell<m<n\le 5).
$$
We also need the 15~planes
$$
F_{k\ell}:=\{x_k=x_\ell=0\}\qquad (0\le k<\ell\le 5)
$$
and the 6~hyperplanes
$$
S_k=\{x_k=0\}\qquad (0\le k\le 5).
$$
Note that
$$
S_k\cap N=\bigcup_{\ell=0}^5 F_{k\ell}.
$$

The map~$p_1$ is the blow up in the 15~points~$P_{k\ell m n}$.
The map~$p_2$ is the blow up along the 20~lines~$L_{k\ell m}^{(1)}$,
where $L_{k\ell m}^{(1)}\subset\PPe$~denotes the strict transform
of~$L_{k\ell m}$
under~$p_1$. More generally, we denote by $L_{k\ell m}^{(i)},
F_{k\ell}^{(i)},\ldots$
the strict transform of $L_{k\ell m},F_{k\ell},\ldots$ in~$\PP^{(i)}$.
Note that the lines $L_{k\ell m}^{(1)}$ are disjoint.
Finally, the map~$p_3$ is the blow up along
the 15~disjoint \lq\lq planes\rq\rq~$F_{k\ell}^{(2)}$. Note that
$F_{k\ell}^{(2)}\cong \hPP_2(6)$, a plane blown up in 6~points.

We define $N^{(i)}$ to be the strict transform of the Nieto quintic
in~$\PP^{(i)}$ and we set $q_i:=p_i|_{N^{(i)}}$.
Note that $q_3$~is an isomorphism.
It was shown in \cite{BN} that the only singularities
of~$\Nd\cong\Nt$ are the 10~Segre points.
The Calabi-Yau model~$Y$ of~$N$ is a small resolution of~$N$,
i.e.\ the map~$q\colon Y\to \Nd$ replaces 10~points by 10~lines.

We now describe the 2:1-coverings. The map
$\pi_i\colon \tPP^{(i)}\to \PP^{(i)}$ is the double cover with branch
locus~$D_i$,
where
$$
D_0:=\sum_{k=0}^5 S_k\subseteq \PPn
$$
and $D_i$~is the odd part of the pull back of~$D_0$
in~$\PP^{(i)}$. Explicitely,
denoting the exceptional divisor of~$p_i$ in~$\PP^{(i)}$ by $E_i$
and its strict transform in~$\PP^{(j)}$ ($j>i$) by $E_i^{(j)}$,
we have
\begin{itemize}
\item $p_1^*(D_0)=\sum S_k^{(1)}+4E_1$,
whence $D_1=\sum S_k^{(1)}$;
\item $p_2^*(D_1)=\sum S_k^{(2)}+3E_2$,
whence $D_2=\sum S_k^{(2)}+E_2$;
\item $p_3^*(D_2)=\sum S_k^{(3)}+E_2^{(3)}+2E_3$,
whence $D_3=\sum S_k^{(3)}+E_2^{(3)}$.
\end{itemize}

\begin{rmk}\label{RX}\upshape Note that $D_3$~and~$\tPPd$ are singular,
since $S_k^{(3)}\cap E_2^{(3)}\neq\emptyset$,
even though Barth and Nieto claim they are smooth \cite[p.~216]{BN}.
This problem disappears after restricting to~$\Nd$, since
$S_k^{(3)}\cap \Nd=\emptyset$.
On the other hand, $S_k^{(2)}\cap \Nt\neq\emptyset$,
so the double cover $\nu_3\colon \tNd\to \Nd\cong\Nt$
is different from $\nu_2\colon \tNt\to\Nt$.
\end{rmk}

Finally, we let $\tN^{(i)}:=N^{(i)}\times_{\PP^{(i)}}\tPP^{(i)}$.
The singular locus of $\tNd$ consists of 20~double points and
$\tq\colon\tY\to \tNd$ is a small resolution.
The branch divisor $D_3\cap \Nd=
E_2^{(3)}\cap \Nd$ of $\nu_3\colon \tNd\to \Nd$ consists
of 20~disjoint quadrics. The map $\hq\colon \tY\to Z$
contracts these quadrics to 20~lines.

\section{The topology of~$Y$ and~$Z$}\label{S2}

In this section we prove the following

\begin{thm}\label{th21} Both $Y$ and~$Z$ have
Hodge numbers $h^{0,0}=h^{3,0}=h^{0,3}=h^{3,3}=1$ and
$h^{1,0}=h^{0,1}=h^{2,0}=h^{0,2}=h^{2,1}=h^{1,2}=
h^{3,1}=h^{1,3}=h^{3,2}=h^{2,3}=0$.
In particular, both manifolds are rigid.
Furthermore,
$h^{1,1}(Y)=h^{2,2}(Y)=50$ whereas $h^{1,1}(Z)=h^{2,2}(Z)=40$.
\end{thm}

\begin{rmk}\upshape This holds for any smooth Calabi-Yau model
of~$Y$ or~$Z$, since birational smooth Calabi-Yau varieties
 have the same Betti numbers \cite{Ba}.
(Note that the Hodge numbers of a smooth Calabi-Yau variety
of dimension~3 are determined by its Betti numbers.)

\end{rmk}

Recall that $Y$ and $Z$ are smooth Calabi-Yau varieties of dimension~3.
This determines the boundary of the Hodge diamond:
$h^{0,0}=h^{3,0}=h^{0,3}=h^{3,3}=1$ and
$h^{1,0}=h^{0,1}=h^{2,0}=h^{0,2}=h^{3,1}=h^{1,3}=h^{3,2}=h^{2,3}=0$.
The remaining Hodge numbers are determined by
$$
a:=h^{2,1}(Y)\qquad\mbox{and}\qquad b:=h^{2,1}(Z)
$$
and the respective Euler characteristics $e(Y)=100$
and $e(Z)=80$. Indeed, we have
$h^{2,1}(Y)=h^{1,2}(Y)=a$,
$h^{2,1}(Z)=h^{1,2}(Z)=b$,
$h^{1,1}(Y)=h^{2,2}(Y)=a+50$ and
$h^{1,1}(Z)=h^{2,2}(Z)=b+40$.

It remains to show that $a=b=0$. This is done in the remainder
of this long section. (A~short roadmap is given
at the end of the next subsection.)
As mentioned in the introduction, we use
the reductions of
 our varieties over suitable finite fields.
We now explain this \lq\lq reduction method\rq\rq\ in detail.

\subsection{The reduction method}

We first deal with~$Y$. Note that $Y$ is defined over~$\ZZ$.
We write $\FF_p:=\ZZ/p\ZZ$; more generally,
we write $\FF_q$ for the finite field with~$q$ elements.

\begin{lem} If $p$ is prime and $p\ge5$, then the reduction
of~$Y$ modulo~$p$ is smooth over~$\FF_p$.
\end{lem}

\begin{pf}{Proof:}{$\boxtimes$} This follows from
an easy calculation (cf.\
 \cite[(3.1), (9.1) and (9.3)]{BN}).

\end{pf}

From now on, $p$~denotes a prime number~$\ge5$.
Let $Y(\FF_p)$ be the set of points of~$Y$ which are
rational over~$\FF_p$.
By Grothendieck's version of the Lefschetz fixed point formula, we have
$$
\#Y(\FF_p)=\sum_{j=0}^6(-1)^j\Tr(\Frob_p^*; \Het^j(Y)),
$$
where~$\Frob_p\colon Y\to Y$ is the Frobenius morphism \cite[p.~454]{H}.
Using the abbreviation
$$
t_j:=\Tr(\Frob_p^*; \Het^j(Y))
$$
we have $t_0=1$, $t_6=p^3$ and $t_1=t_5=0$, whence
$$
\#Y(\FF_p)=1+t_2-t_3+t_4+p^3.
$$
The Weil conjectures assert that the eigenvalues of~$\Frob_p$
on~$\Het^j(Y)$ are algebraic integers, which do not depend
on~$\ell$ and which have (archimedian) absolute value~$p^{j/2}$.
This implies that the $t_j$ are ordinary integers,
independent of~$\ell$, with absolute value
$|t_j|\le b_j(Y)\cdot p^{j/2}$.
We use this only for $j=3$:
$$
|t_3|\le (2a+2)p^{3/2}.
$$
For $j=2$ we use the following stronger result, the proof of which
we postpone to the end of this section.

\begin{pro}\label{p24} If $p\ge5$ then
all eigenvalues of~$\Frob_p$ on~$\Het^2(Y)$ are equal to~$p$.
\end{pro}
Hence $t_2=b_2(Y)p=(a+50)p$ and $t_4=(a+50)p^2$, by Poincar\'{e} duality.
Putting everything together, we get
$$
|1+(a+50)(p+p^2)+p^3-\#Y(\FF_p)|\le (2a+2)p^{3/2}.
$$
For $p=13$ we will find $\#Y(\FF_p)=11260$, whence
$182a+38\le 94(a+1)$ implying $a=0$.

We use the same method to show that $b=0$.
In fact, we work with~$\tY$ rather than~$Z$.
The reason is that
it is not clear to us if
the map~$\hq\colon\tY\to Z$ is defined over~$\ZZ$,
since it was obtained in \cite{BN} using Mori theory.
On the other hand,
$\tY$ is clearly defined over~$\ZZ$.
Furthermore, $h^{i,j}(\tY)=h^{i,j}(Z)$
for $(i,j)\neq (1,1)$ or $(2,2)$, whereas
$h^{1,1}(\tY)=h^{2,2}(\tY)=b+60$.

\begin{lem} If $p$ is prime and $p\ge5$, then the reduction
of~$\tY$ modulo~$p$ is smooth over~$\FF_p$.
\end{lem}

\begin{pf}{Proof:}{$\boxtimes$}
Just note that
the proof of \cite[(10.2)]{BN} still works.
\end{pf}

\noindent
We also need the following proposition, which we will
also prove at the end of this section.

\begin{pro}\label{p26} If $p\equiv 1 \bmod 4$ then
all eigenvalues of~$\Frob_p$ on~$\Het^2(\tY)$ are equal to~$p$.
\end{pro}
So for $p\equiv 1\bmod 4$ (which implies $p\ge5$) we have
$$
|1+(b+60)(p+p^2)+p^3-\#\tY(\FF_p)|\le (2b+2)p^{3/2}.
$$
For $p=13$ we will see that $\#\tY(\FF_p)=13080$, whence
$182b+38\le 94(b+1)$ implying $b=0$.

\

Summarizing, we have to determine the number of $\FF_{13}$-rational
points on $Y$ and~$\tY$, and we have to prove Propositions \ref{p24}
and~\ref{p26}.
This is done in \S \ref{ss23} and~\S \ref{ss24}, respectively.
The next subsection is a preparation for \S\ref{ss23}.

\subsection{Counting points on Cayley cubics}

In order to count the points of $Y(\FF_p)$
we will need to understand the structure of the
exceptional divisor in~$Y$ lying over a point~$P_{k\ell m n}$.
We now collect the necessary information.
Since $Y\to\Nd\cong\Nt$ is an isomorphism outside
the Segre points, we can consider~$\Nt$
instead of~$Y$: the fibres over~$P_{k\ell m n}$ are the same.
One easily checks (see \cite[(9.1)]{BN} and the proof of \ref{L28} below)
that the fibre of $\Ne\to N$ lying over~$P_{0123}=(0:0:0:0:-1:1)$
is the Cayley cubic
$$
\Ce\colon \qquad \frac1{y_0}+\frac1{y_1}+\frac1{y_2}+\frac1{y_3}
$$
in~$\PP_3$. This $\PP_3$ is of course the component of~$E_1\subset\PPe$
lying over~$P_{0123}$ and $\Ce=\PP_3\cap \Ne$.

Analogously, in order to determine $\#\tY(\FF_p)$
we need information about the fibre of~$\tY\to N$
lying over~$P_{k\ell m n}$. Since $\tY\to \tNd$
is an isomorphism  outside the Segre points,
we may replace $\tY$ by~$\tNd$.
Now recall that the double cover $\tNd\to \Nd\cong\Nt$
is branched along the exceptional divisor~$E_2\cap \Nt$
of the blow up $\Nt\to\Ne$ along the lines $L_{k\ell m n}^{(1)}$.
Hence $\tCt$ is the double cover of~$\Ct$ branched along
$E_2\cap \Ct$, which consists of the four exceptional
lines of~$\Ct\to\Ce$.

\begin{lem}\label{L28} Let $p\equiv 1\bmod 4$.
Then $\#\tCd(\FF_p)=p^2+8p+1$.
\end{lem}

\begin{pf}{Proof:}{$\boxtimes$}
First we determine the blow up of $N$ in the point $P_{0123}=(0:0:0:0:-1:1)$.
We blow up the $\AAA_5\subset\PP_5$ where $x_5=1$ in $P=(0,0,0,0,-1)$,
this is the subvariety of $\AAA_5\times\PP_4$ defined by
$$
x_iy_j-x_jy_i=0,\qquad (x_4+1)y_j-x_jy_4=0,\qquad 0\leq i<j\leq 3.
$$
The inverse image of $N\cap \AAA_5$ in the open set $\AAA_5\times\AAA_4$
defined by $y_4=1$ is given by the equations $x_j=(x_4+1)y_j$, so we are left
with the variables $x_4,\,y_0,\ldots,\,y_3$, the equation
$$
(x_4+1)(y_0+y_1+y_2+y_3+1)=0
$$
(obtained from $\sum x_i=0$) and the equation
$$
(x_4+1)^3\left(
y_0y_1y_2y_3(x_4+1)^2+x_4(y_0y_1y_2+y_0y_1y_3+y_0y_2y_3+y_1y_2y_3)\right)=0
$$
(obtained from $\sigma_5=x_0x_1x_2x_3x_4+\cdots+x_1x_2x_3x_4x_5=0$).
The fibre over~$P$ is defined by $x_4=-1$ and
we see that the fibre over $P_{0123}$ in the blow up $N^{(1)}\rightarrow N$
is the Cayley cubic $C^{(1)}$ in a $\PP^3$
(note $y_4=-\sum y_i)$ with coordinates $y_0,\ldots,y_3$.

The double cover $\tN$ of $N$ ramifies over the zero locus of
$x_0x_1x_2x_3x_4x_5$.
One can thus construct it naturally in a weighted projective space.
Equivalently, consider the image of $N$ in $\PP^k$ under the third Veronese
map,
so the coordinates on $\PP^k$ are the $u_{ijk}$ with $u_{ijk}=x_ix_jx_k$.
Then $\tN$ is subvariety of $\PP^{k+1}$, with coordinates $u_{ijk}$ and $v$,
defined by the equations of the image of $N$ in $\PP^k$ and the equation
$v^2=u_{012}u_{345}$.
The inverse image of $N\cap \AAA_5$ in $\tN$ is therefore isomorphic
to the subvariety of $\AAA_6$, with coordinates $x_0,\ldots,\,x_4$ and $v$,
defined by the equations of $N\cap \AAA_5$ and the equation
$v^2=x_0x_1x_2x_3x_4$.

To obtain $\tN^{(1)}$, locally near $P_{0123}$,
we consider the inverse image of $N\cap \AAA_5$
in the open set of the blow up where $y_4=1$:
$$
v^2=(x_4+1)^4x_4y_0y_1y_2y_3.
$$
As explained in section \ref{sec1},
we need to ramify only over the odd part of this divisor,
so locally $\tN^{(1)}$ is defined by
$$
w^2=x_4y_0y_1y_2y_3
$$
and the equations defining the strict transform of $N\cap \AAA_5$.
Restricted to the special fibre $x_4=-1$ this gives the equation
$w^2=-y_0y_1y_2y_3$. In case $-1$ is a square in the field under consideration
(for example $\FF_p$ with $p\equiv 1$ mod $4$),
we can change the variable to obtain the equation $w_1^2=y_0y_1y_2y_3$.

The equation $y_i=0$ defines the divisor $S_i^{(1)}$. The intersection
$S_i^{(1)}\cap C^{(1)}$ consists of three lines in $C^{(1)}$.
In this way we get the $6$ lines $S^{(1)}_i\cap S^{(1)}_j\cap C^{(1)}$
in the Cayley cubic which connect the nodes. In particular, the cover
$\tCe\rightarrow \Ce$ ramifies over these $6$ lines and these lines are
singular on $\tCe$.

The map $\Ct\rightarrow \Ce$ is the blow up of the Cayley cubic
in the $4$ nodes.
We denote the exceptional fibres, $(-2)$-curves, by $R_i^{(2)}$.
The double cover $\tCt\rightarrow \Ct$ ramifies along the $6$
lines and the $4$ $R_i^{(2)}$'s. The map $N^{(3)}\rightarrow N^{(2)}$ is an
isomorphism
(since the $F_{kl}^{(2)}$ are smooth surfaces in the smooth part
of the threefold $N^{(2)}$). But since $F_{kl}^{(1)}=S_k^{(1)}\cap S_l^{(1)}$
intersects the Cayley cubic in a line, the map $\tCd\rightarrow \tCt$
is the normalization of $\tCt$. Therefore the double cover
$$
\tCd\longrightarrow \Cd=\Ct
$$
ramifies only over the four nodal curves $R_i$.
It is easy to check that $h^{1,0}(\tCd)=h^{2,0}(\tCd)=0$
and $h^{1,1}(\tCd)=8$.

We recall the explicit parametrizations for the Cayley cubic and its double
cover,
ramified only in the nodes. The cubic $\Ce$ is obtained by blowing up
a~$\PP_2$ in $6$~points
which are the intersection points of $4$ general lines, next one blows down
the strict transforms of the $4$ lines.
For the plane we will take
$$
\PP_2:=\{x\in\PP_3\colon x_0+x_1+x_2+x_3=0\},
$$
the $4$ lines are defined by $x_i=0$ so
the 6~points $x_i=x_j=0$ ($0\le i<j\le 3$).
The linear system of cubics through these 6~points
is generated by $x_0x_1x_2,x_0x_1x_3,x_0x_2x_3$,
and~$x_1x_2x_3$.
Denoting temporarily
the coordinates of the target space~$\PP_3$ by $(y_0:y_1:y_2:y_3)$,
we see that
the associated rational map ${\PP_2 \dashrightarrow \PP_3}$,
$(x_0:x_1:x_2:x_3)\mapsto
x_0x_1x_2x_3(\frac1{x_0}:\frac1{x_1}:\frac1{x_2}:\frac1{x_3})=:(y_0:y_1:y_2:
y_3)$,
maps $\PP_2$ to the Cayley cubic $\sum \frac1{y_i}=0$.
It blows up the six~points, the exceptional divisors
being the edges of the tetrahedron spanned by
the four singular points of~$\Ce$,
and it contracts the four~lines~$x_i=0$
to these singular points.
In other words, the map $\PP_2\dashrightarrow \Ce$
induces an isomorphism
$$
\hPP_2(6)\stackrel{\sim}{\longrightarrow}\Ct,
$$
where $\hPP_2(6)$ is the blow up of $\PP_2$ in the $6$ points.
This also works over the finite field~$\FF_p$,
hence $\#\Ct(\FF_p)=1+7p+p^2$.

An alternative description of~$\tCd$
would be the double cover of~$\PP_2$ branched
along the four lines~$x_i=0$ with the six quadratic
singular points lying over the intersection points resolved
by a $(-2)$-curve.

Since each of the $6$ points lies on exactly two lines,
the inverse image in $\tCd$ of the exceptional divisor is an
irreducible rational curve which maps $2:1$ onto the exceptional divisor in
$\Cd$.
Note that the line $x_3+x_4=0$ is parametrized by $(s:-s:t:-t)$ and
splits in this double cover since it meets the ramification locus only
in the points $(1:-1:0:0)$ and $(0:0:1:-1)$ and each point has multiplicity
$2$.

Consider the morphism $f\colon\Cd\rightarrow \PP_1$ defined by the pencil of
lines on the point $(1:-1:0:0)$
(so $f$~is obtained from the rational map
$\PP_3\dashrightarrow \PP_1$, $(x_0:x_1:x_2:x_3)\mapsto (x_2:x_3)$).
If the line in the pencil~$f$ meets the ramification locus in two other
distinct points,
besides $(1:-1:0:0)$, the corresponding fibre of
the composite map~$\tf\colon \tCd\to\Cd\to\PP_1$ is smooth.
There are only three exceptions: the lines $x_2=0$, $x_3=0$ and $x_3+x_4=0$.

The fibre of~$\tf$ corresponding to the line~$x_3=0$
consists of the strict tranform of this line as well the two $(-2)$-curves
corresponding to the points $(1:0:-1:0)$ and $(0:1:-1:0)$. Thus $\tf^{-1}(1:0)$
is a tree of three $\PP^1$'s, and same holds for $\tf^{-1}(0:1)$.
We already observed that the line $x_3+x_4=0$
splits in the covering. Besides these two components, $f^{-1}(1:-1)$ also
contains
the $(-2)$-curve over the point $(0:0:1:-1)$ and this fibre is thus also a
tree of three $\PP^1$'s.

If the two components over the line $x_3+x_4=0$ are rational over $\FF_p$,
then the number of points of $\tCd$ over $\PP^1$ minus the three
exceptional points is $(p-2)(p+1)$,
and each exceptional point contributes a tree of three
$\PP^1$'s hence $3(p+1)-2=3p+1$ points. Thus for such fields
$$
\# \tCd(\FF_p)=(p^2-p-2)+3(3p+1)=p^2+8p+1.
$$
From this one easily obtains that $\dim H^2(\tCd)=8$.
(On $\tCd$ we have $7$~obviously independent divisor classes,
the pull-back of a hyperplane section of $\Cd$ and the $6$ nodal curves.
Each of the components over the line $x_3+x_4=0$
can serve as an \lq\lq eighth divisor\rq\rq.
To see this, one checks that the $8\times 8$ intersection matrix
has non-zero determinant.)

To see for which primes the two components over $x_3+x_4=0$ are rational over
$\FF_p$, we pull-back the defining equation for the double cover
$w^2=-y_0y_1y_2y_3$ (and $y_4=-(y_0+\ldots+y_3)$) to $\PP_2$, and one finds
that
that over $(1:-1:t:-t)$ the double cover is given by $w^2=-t^2$,
hence the two components are rational iff $-1$ is a square in the field iff
$p\equiv 1$ mod $4$.

\end{pf}

\begin{rmk}\upshape
For completeness sake we give the well-known identification
of $\tCd$ with a $\PP^2$ blown up in $7$ points. The standard
Cremona transformation on $P_2$,
$(t_0:t_1:t_2)\mapsto (1/t_0:1/t_1:1/t_2)$,
is well defined on the blow up of $P^2$ in the three `basis' points and has
4~fixed points: $(1:1:1), (1:1:-1), (1:-1:1), (-1:1:1)$.
Blowing up these $7$ points gives a surface with $1+p+p^2+7p=1+8p+p^2$ points.
The smooth quotient surface has 4 (-2)-curves (the images of the exceptional
divisors over the fixed points) and blowing these down you get the Cayley
cubic $\Ce$.
Therefore $\tCd$ is the isomorphic with the blowup of $\PP^2$ in $7$ points.
An explict formula for the rational 2:1 map from $\PP^2$ to $\Ce$
is given by the following polynomials, which are invariant (modulo the
common factor $(t_0t_1t_2)^2$) under the Cremona transformation
$$
x_0=2t_0t_1t_2,\quad
x_1=t_0(t_1^2+t_2^2),\quad
x_2=t_1(t_0^2+t_2^2),\quad
x_3=t_2(t_0^2+t_1^2)
$$
and they satisfy the cubic relation
$2x_1x_2x_3-x_0(x_1^2+x_2^2+x_3^2-x_0^2)$.
The fixed points map to
$(1:1:1:1)$, $(1:1:-1:-1)$, $(1:-1:1:-1)$ and $(1:-1:-1:1)$.
Using the linear combinations
$x_0+x_1+x_2+x_3$,  $x_0-x_1+x_2-x_3$, $x_0+x_1-x_2-x_3$ and
$x_0-x_1-x_2+x_3$, one obtains functions which
satisfy the equation defining $\Ce$.
\end{rmk}

\begin{cor}\label{C29}
All eigenvalues of~$\Frob_p$ acting on~$\Het^2(\tCt)$
are equal to~$p$
if $p\equiv 1\bmod 4$.
\end{cor}

\begin{pf}{Proof:}{$\boxtimes$} The eigenvalues
$\alpha_1,\ldots,\alpha_8$
of~$\Frob_p$ acting on~$\Het^2(\tCt)$ are
complex numbers of absolute value~$p$.
Since $\alpha_1+\cdots+\alpha_8=8p$ by the previous lemma,
we find that $\alpha_i=p$ for $i=1,\ldots,8$.
\end{pf}

\begin{rmk}\upshape For $p\equiv 3\bmod 4$ one finds
$\#\tCt(\FF_p)=p^2+6p+1$.
Thus the eigenvalues of~$\Frob_p$ acting on~$\Het^2(\tCt)$
are $p$~with multiplicity~7 and $-p$~with multiplicity~1.
\end{rmk}

\begin{lem}\label{L212} If $p\equiv 1\bmod 4$ then
the contribution over the unramified
part of $F_{k\ell}$ to $\#\tY(\FF_p)$ is
$p^2-2p+3$.
\end{lem}

\begin{pf}{Proof:}{$\boxtimes$}
We claim that the equation of
the unramified part of the double cover of~$F_{k\ell}$
is $y^2=-x_0x_1x_2x_3$ as in Lemma~\ref{L28}.
In order to see this, we blow up~$\PP^{(0)}$ along $F_{k\ell}$.
The resulting exceptional divisor is $\PP_1\times F_{k\ell}$
and one easily checks that the strict transform of~$N$
intersects this divisor in~$\{(1:-1)\}\times F_{k\ell}$
with respect to the coordinates on the fibre~$\PP_1$
given by $x_k$~and~$x_\ell$. On $\{(1:-1)\}\times F_{k\ell}$
we have $x_k=-x_\ell$, so
restricting
$x_0\cdots x_5$
and dividing by the square~$x_k^2=x_\ell^2$ yields
the desired $-x_0\cdots x_5/x_kx_\ell$.

Now we can apply Lemma~\ref{L28}.
Indeed, the ramification locus, which consists
of 6~intersecting lines, contains $10p-2$ points over~$\FF_p$,
the proper double cover~$\widetilde{F}_{k\ell}\cong\tCt$ of~$F_{k\ell}$
contains $p^2+8p+1$ points, leaving $(p^2+8p+1)-(10p-2)$
points for the unramified part.
\end{pf}

\begin{rmk}\upshape The corresponding result for $p\equiv 3\bmod 4$, $p\neq3$,
is $p^2-4p+3$.
\end{rmk}

\subsection{Counting points on~$Y$ and~$\tY$}\label{ss23}

We will now show how to count the points
of~$Y$ and~$\tY$ over~$\FF_p$.
Recall that $D_0=\{x_0\cdots x_5=0\}$. We set
$$
U:=N\setminus D_0.
$$
It is easy to count the points of $U(\FF_p)$ on a computer.
For $p=13$ one finds $\#U(\FF_{13})=2140$.
Together with the following proposition
this yields $\#Y(\FF_{13})=11260$.

\begin{pro} For $p\ge5$ we have
$\#Y(\FF_p)=\#U(\FF_p)+50p^2+50p+20$.
\end{pro}

\begin{pf}{Proof:}{$\boxtimes$}
Recall that $N\cap D_0$ is the union of the 15~planes~$F_{k\ell}$.
The map $Y\to N$ replaces the 10~Segre points by 10~lines,
the 20~lines~$L_{k\ell m}$ by quadrics and the
15~points~$P_{k\ell m n}$ by smooth cubics~$\Ct\cong\hPP_2(6)$.
The Segre points lie in~$U$, so here we get $10p$ new points over~$\FF_p$.
Over $R^0:=\cup F_{k\ell}\setminus \cup L_{k\ell m}$ nothing happens.
Now $R^0$ consists of 15~copies of $\PP_2\setminus \mbox{(4~lines)}$,
hence
$$
\#R^0(\FF_p)=15[(p^2+p+1)-4(p+1)+6]=15(p^2-3p+3).
$$
Similarly, $\cL^0:=\cup L_{k\ell m}\setminus \cup P_{k\ell m n}$
consists of 20~copies of $\PP_1\setminus\mbox{(3~points)}$.
Now $Y\to N$ replaces $\cL^0$ by $\cL^0\times \PP_1$ and
its contribution to~$\#Y(\FF_p)$ is
$$
\#(\cL^0\times\PP_1)(\FF_p)=20(p-2)(p+1).
$$
Finally, we have a copy of~$\hPP_2(6)$ over each of the
15~points~$P_{k\ell m n}$, contributing
$$
15(p^2+7p+1).
$$
Hence
\begin{eqnarray*}
\#Y(\FF_p)&=&\#U(\FF_p)+10p+15(p^2-3p+3)\\
&&\qquad +20(p-2)(p+1)+15(p^2+7p+1)\\
&=& \#U(\FF_p)+50p^2+50p+20.
\end{eqnarray*}
\end{pf}

Now we count the points of~$\tY(\FF_p)$.
Let $\tU\subset \tY$ be the inverse image of~$U\subset Y$.
One can again easily count the points in~$\tU(\FF_p)$
using
$$
\#\tU(\FF_p)=2\cdot\#\{x\in U(\FF_p)|\mbox{$x_0\cdots x_5$ is a square
in~$\FF_p$}\}
$$
and a computer. Together with the following proposition
this yields the desired $\#\tY(\FF_{13})=2\times 1720+9640
=13080$.

\begin{pro} If $p\equiv 1\bmod 4$ then
$\#\tY(\FF_p)=\#\tU(\FF_p)+50p^2+90p+20$.
\end{pro}

\begin{pf}{Proof:}{$\boxtimes$}
Note that the 20~Segre points in~$\tNd$ lie in~$\tNd(\FF_p)$,
since $1^3\cdot (-1)^3$ is a square in~$\FF_p$. Hence
$\#\tY(\FF_p)=\#\tNd(\FF_p)+20p$
and we have to show that
$\#\tNd(\FF_p)=\#\tU(\FF_p)+50p^2+70p+20$.
Now recall that the branch locus of $\tNd\to\Nd=\Nt$ is~$E_2^{(2)}$,
the union of 20~disjoint quadrics.
The intersection of this branch locus with the resolved
Cayley cubic~$\Ct$ in~$\Nt$ lying over a point~$P_{k\ell mn}$
consists of the 4~exceptional lines of $\Ct\to\Ce$.
We refer to the complement
as the open Cayley cubic.
The contribution of the ramification locus~$E_2^{(2)}$ to~$\#\tNd(\FF_p)$
is $20(p+1)^2$. By Lemma~\ref{L28}
each open Cayley cubic contributes $(p^2+8p+1)-4(p+1)=
p^2+4p-3$.
Lemma~\ref{L212} implies that the contribution of
$R^0=\cup F_{k\ell}\setminus \cup L_{k\ell m}$
equals $15(p^2-2p+3)$.
So we find
$
\#\tNd(\FF_p)-\#\tU(\FF_p)=20(p+1)^2+15(p^2+4p-3)+15(p^2-2p+3)
=50p^2+70p+20
$.
\end{pf}

\begin{rmk}\upshape For $p\equiv 3\bmod 4$, $p\neq 3$, one
finds
$\tY(\FF_p)=
\#\tNd(\FF_p)=\#\tU(\FF_p)+20(p+1)^2+15(p^2+2p-3)+15(p^2-4p+3)
=\#\tU(\FF_p)+50p^2+10p+20
$.
\end{rmk}

\subsection{Proofs of Propositions \ref{p24} and~\ref{p26}}\label{ss24}

In the remainder of this section we prove Propositions
\ref{p24} and~\ref{p26}
about the action of Frobenius on $\Het^2(Y)$ and $\Het^2(\tY)$.
We say that a variety over~$\ZZ$
satisfies condition~$(*)_p$ for a prime~$p$
 if all eigenvalues of~$\Frob_p$
acting on the second \'etale cohomology are equal to~$p$.
If $(*)_p$ holds for every good prime~$p$, then
we say that the variety satisfies condition~$(*)$;
if $(*)_p$ holds for every good prime~$p\equiv 1\bmod 4$,
then we say that the variety satisfies condition~$(**)$.
So we have to show that $Y$ and~$\tY$ satisfy $(*)$ and~$(**)$,
respectively. (Note, however, that for the proofs of $a=b=0$ above,
we only used the fact that $Y$ and~$\tY$ satisfy condition~$(*)_p$
for $p=13$.)

We start with~$Y$.
Recall the sequence of maps
$Y\stackrel{q}{\to} \Nd=\Nt
\stackrel{q_2}{\to}\Ne \stackrel{q_1}{\to} \Nn=N\subset\PP_4$.
Since $q^*\colon \Het^2(\Nd)\to \Het^2(Y)$ is an isomorphism,
it suffices to check that condition~$(*)$ holds for~$\Nd=\Nt$.
So the following lemma finishes the proof of Proposition~\ref{p24}.

\begin{lem}\label{L17} Condition~$(*)$ holds for $N=\Nn,\Ne$ and $\Nt=\Nd$.
\end{lem}

\begin{pf}{Proof:}{$\boxtimes$} The condition holds for~$\Nn=N$,
since $\Het^2(N)=\Het^2(\PP_4)$ by the Lefschetz hyperplane
theorem \cite[Corollary~I.9.4]{FK}.
In order to lift this result to~$\Ne$, we use the spectral sequence
$$
E_2^{j,2-j}=\Het^j(\Nn,R^{2-j}q_{1,*}(\QQ_\ell))\implies \Het^2(\Ne)
$$
associated to $q_1\colon \Ne\to \Nn$.
Since condition~$(*)$ is stable under extensions
of Galois modules,
it suffices to check the analogous condition on the graded
parts~$E_2^{j,2-j}$ of~$\Het^2(\Ne)$.
Now we already know that Frobenius acts correctly on $E_2^{2,0}=\Het^2(\Nn)$.
For $j=1$ we find
$E_2^{1,1}=\Het^1(\cP,R^1q_{1,*}(\QQ_\ell))=0$,
since the support $\cP:=\{P_{k\ell m n}\}$ of~$R^1q_{1,*}(\QQ_\ell)$
has dimension~0.
Finally, we have to consider $E_2^{0,2}=\Het^0(\cP,R^2q_{1,*}(\QQ_\ell))
= \Het^2(\Ce)^{\oplus 15}$.
So it suffices to note that $\Ce$ satisfies condition~$(*)$.
Indeed, $H^2(\Ce)$~is a Galois submodule of~$H^2(\Ct)$
and the latter space is generated by a line
and 6~exceptional divisors, all of which are
obviously defined over~$\FF_p$.

The step from $\Ne$ to $\Nt$ is similar.
Recall that $\Nt\to\Ne$ contracts 20~quadrics to lines.
Hence we have a spectral sequence
$$
E_2^{j,2-j}=\Het^j(\Ne,R^{2-j}q_{2,*}(\QQ_\ell))\implies \Het^2(\Nt)
$$
with $E_2^{2,0}=\Het^2(\Ne)$, $E_2^{1,1}=0$ and $E_2^{0,2}=\oplus
\Het^2(L_{k\ell m}^{(1)})
=\QQ_\ell(-1)^{\oplus 20}$. This proves the assertion for~$\Nt$.
\end{pf}

Finally, we prove Proposition~\ref{p26}, i.e.\ we check
condition~$(**)$ for~$\tY$.
Since the Leray spectral sequence of $\tq\colon \tY\to \tNd$ induces
an exact sequence
$$
0\to \Het^2(\tNd)\to \Het^2(\tY) \to  H^2(\PP_1)^{\oplus 20},
$$
where the $\PP_1$'s the fibres over the Segre points,
is suffices to show that $\tNd$~satisfies~$(**)$.
(Note that the $\PP_1$'s are defined over~$\FF_p$.)
We use a Lefschetz argument relating
$\Het^2(\tNd)$ to $\Het^2(H)$ for a suitable
divisor $H\subset \tNd$.
The complement~$\tNd\setminus H$ will contain
the 20~Segre points, so we have to modify the
hyperplane theorem to deal with them.
This is done in the following lemma, which we will
apply to~$X=\tNd$.

\begin{lem}\label{L217} Let $X$ be a projective threefold.
Let $H\subset X$ be a divisor such that
$X\setminus H$~is affine and $\Sigma:=\Sing(X\setminus H)$ consists
of at most finitely many $A_1$-singularities.
Let $\hX\to X$ be the blow up $X$ in~$\Sigma$
and let $Q$~be the union of the exceptional
quadrics. Consider the natural maps
$\alpha\colon H^2(\hX)\to H^2(H)$ and $\beta\colon H^2(\hX)\to H^2(Q)$.
Then $\beta$~is injective on the kernel of~$\alpha$.
\end{lem}

\begin{pf}{Proof:}{$\boxtimes$} Note that
$H^2(H\cup Q)=H^2(H)\oplus H^2(Q)$ since
$H\cap Q=\emptyset$. It follows that $H^2(H\cup Q,H)\cong H^2(Q)$
and we have a commutative diagram
$$
\begin{CD}
@. H^2(\hX,H\cup Q)\\
@. @VVV\\
H^1(H) @>>> H^2(\hX,H) @>>> H^2(\hX) @>{\alpha}>> H^2(H)\\
@. @VVV @VV{\beta}V\\
@. H^2(H\cup Q,H) @>{\sim}>> H^2(Q)\\
@. @VVV\\
@. H^3(\hX,H\cup Q)
\end{CD}
$$
In order to prove the claim it is sufficient to show that
$H^2(\hX,H\cup Q)=0$. By Alexander duality
$$
H^2(\hX,H\cup Q) = H_4(\hX\setminus (H\cup Q)) \cong H_4(W\setminus\Sigma),
$$
where $W:=X\setminus H$.  Since $W$~is affine of dimension~3
we know that $H_k(W)=0$ for $k\ge4$. The exact homology sequence
for the pair $(W,W\setminus\Sigma)$ gives
$$
 H_5(W,W\setminus\Sigma) \to H_4(W\setminus \Sigma) \to H_4(W).
$$
Since $H_4(W)=0$ it is enough to prove that $H_5(W,W\setminus\Sigma)=0$.
Let us first assume that $\Sigma=\{p\}$. Then by an excision argument
$H_5(W,W\setminus\Sigma)\cong H_5(U,U\setminus p)$ where $U$ is a
suitable contractible neighbourhood of the point~$p$.
It remains to show that $H_5(U,U\setminus p)=0$.
This is a local problem. Let $Q'\subset\CC^4$
be the quadric cone and let $K=S_\epsilon \cap Q'$
be the link.
Since $U$~is contractible, we find $H_5(U,U\setminus p)\cong H_4(K)$.
The manifold~$K$ is orientable and 1-connected
(see \cite[p.~76]{Di}).
Hence $H_4(K)\cong H_1(K)\cong 0$
where the first isomorphism is by Poincar\'{e} duality.
The general case where $\Sigma$~consists of several points can be done by
induction by removing one point at a time.
\end{pf}

\begin{rmk}\upshape Note that the proof works
for arbitrary isolated hypersurface singularities.
\end{rmk}

As mentioned before, we will apply the lemma to $X=\tNd$.
For $H$ we will take the pull back of
$$
h:=E_1^{(3)}+E_2^{(3)}+E_3.
$$

\begin{lem} The pair $(\tNd,H)$ satisfies the conditions of Lemma~\ref{L217}.
\end{lem}

\begin{pf}{Proof:}{$\boxtimes$}
It suffices to show that
$a_0\sum S_k^{(3)}+a_1E_1^{(3)}+a_2E_2^{(3)}+a_3E_3$ is very ample
on~$\PPd$ for suitable integers $a_0,a_1,a_2,a_3>0$.
Indeed, its restriction to~$\Nd$ has the same support
as~$h$, since $S_k^{(3)}\cap \Nd=\emptyset$.
Hence its pull back~$H$ supports
an ample divisor on~$\tNd$, since $\tNd\to\tN$ is finite.
This implies that the complement $\tNd\setminus H$ is affine.
We also know that $\tNd\setminus H$ is non-singular outside
the Segre double points.

It remains to prove the claim about
$a_0\sum S_k^{(3)}+a_1E_1^{(3)}+a_2E_2^{(3)}+a_3E_3$.
Let $H_0:=\sum S_k$ in~$\PPn$. The divisor
$H_1:=p_1^*(n_1H_0)-E_1$ is very ample on~$\PPe$
for $n_1\gg0$. In fact it suffices that $\cI_\cP(6n_1-1)$
be globally generated, where $\cP=\{P_{k\ell m n}\}$
(see \cite[proof of Proposition~II.7.10~(b)]{H}).
Similarly, $H_2:=p_2^*(n_2H_1)-E_2$ is very ample on~$\PPt$
for $n_2\gg n_1\gg 0$ and
$H_3:=p_3^*(n_3H_2)-E_3$ is very ample on~$\PPd$
for $n_3\gg n_2\gg n_1\gg 0$.
Using $p_1^*(\sum S_k)=\sum S_k^{(1)}+4E_1$,
one finds that $H_1=n_1\sum S_k^{(1)}+(4n_1-1)E_1$.
Continuing in this fashion for $H_2$ and~$H_3$,
one finds $H_3=a_0\sum S_k^{(3)}+a_1E_1^{(3)}+a_2E_2^{(3)}+a_3E_3$
with $a_0=n_3n_2n_1$, $a_1=n_3n_2(4n_1-1)$,
$a_2=n_3(3n_2n_1-1)$ and $a_3=2n_3n_2n_1-1$.
So $H_3$ is indeed very ample for suitable $a_0,a_1,a_2,a_3>0$.
\end{pf}

\begin{lem} Let $Q$ be an exceptional quadric obtained
by blowing up~$\tNd$ in a Segre node. Then
condition~$(**)$ holds for~$Q$.
\end{lem}

\begin{pf}{Proof:}{$\boxtimes$} We may assume that the node is
$(1:1:1:-1:-1:-1)$.
The tangent cone,
which is a cone over~$Q$, is given by the equation
$$
x_0+\cdots+x_5=x_0^2+x_1^2+x_2^2-x_3^2-x_4^2-x_5^2=0
$$
in~$\PP_5$ (see \cite[p.~181]{BN}). So the rulings
of~$Q$ are defined over~$\FF_p$.
\end{pf}

\begin{pro} Condition~$(**)$ holds for~$\tNd$.
\end{pro}

\begin{pf}{Proof:}{$\boxtimes$}
Let $X=\tNd$ and let $\hX$ be the blow up of~$X$
in the 20~Segre nodes. The Leray spectral sequence of
$\hX\to X$ induces an injection $H^2(X)\into H^2(\hX)$,
so it suffices to prove condition~$(**)$ for~$\hX$.
By Lemma~\ref{L217},
the comparison theorem between \'{e}tale and complex cohomology,
and the lemma above
it would even suffice to check
condition~$(**)$ for~$H$.
However, we will see in a moment that
condition~$(**)$ probably does not hold for~$H$.
This is not a serious problem.
Indeed, we just restrict our attention
to eigenvalues {\em of absolute value~$p$}
in the conditions $(*)_p$,
$(*)$ and~$(**)$. Call the (weaker) modified
conditions $(*)_p'$, $(*)'$ and~$(**)'$, respectively.
For smooth proper varieties
the modified conditions are equivalent to the original ones.
Furthermore, the modified conditions are again
stable under extensions of Galois modules.
So it suffices to prove that $(**)'$ holds for~$H$.

Note that $h=E_1^{(3)}\cup E_2^{(3)}\cup E_3$ is a divisor
with global normal crossings, i.e.\ its components
are smooth and meet each other transversally.
The same holds for the pull back~$H$.
So we can apply the Mayer-Vietoris spectral
sequence (cf.\ \cite[p.~103]{M}).
This spectral
sequence abuts to $H^{p+q}(H)$,
degenerates at~$E_2$
and has $E_1^{p,q}=H^q(X[p])$,
where $H[i]$~is the disjoint union of strata
of codimension~$i$, i.e.\ $H[0]$~consists of the triple points,
$H[1]$~is the disjoint union of the intersection curves
and $H[0]$~is the disjoint union of the smooth components
of~$H$.
Now note that the components of~$H[1]$ are rational curves,
hence $E_1^{1,1}=0$.
So in our case we get an exact sequence
$$
H^0(H[1])\stackrel{d_1^{1,0}}\to H^0(H[2]) \to H^2(H)
\to H^2(H[0])\stackrel{d_1^{0,2}}{\to} H^2(H[1]).
$$
This exact sequence is related to the weight filtration
on~$H^2(H)$: one has $W_0(H^2(H))=W_1(H^2(H))=
E_\infty^{2,0}=\Coker(d_1^{1,0})$ and $\Gr_2^W(H^2(H))=
E_\infty^{0,2}=\Ker(d_1^{0,2})$
(cf.\ \cite{M}).
Corollary~\ref{C29} implies that condition~$(**)$
holds for $\Gr_2^W\subseteq H^2(H[0])$.
On the other hand, $\Frob_p={\mathrm {id}}$ on~$W_0$
(cf.\ \cite[\S14]{De2}), so condition~$(**)$ does not
hold on~$W_0$ if $W_0\neq0$!
However, as we remarked in the first paragraph of this proof,
it suffices to note that the weaker condition~$(**)'$
clearly holds for~$W_0$.

\end{pf}

\begin{rmk}\upshape
Alternatively, one can use the strictness of the weight filtration
(cf.\ \cite[\S1]{De2}). Indeed, we want to show that
for $p\equiv 1\bmod 4$
the eigenvalues of $\Frob_p$ on the image~$I$
of $\alpha\colon H^2(\hX)\to H^2(H)$ are all
equal to~$p$. Strictness of the weight filtration
implies that $I\cap W_1=0$, i.e.\
that $I$~injects into~$\Gr_2^W(H^2(H))\subseteq H^2(H[0])$,
where we know condition~$(**)$ by Corollary~\ref{C29}.
\end{rmk}

\section{The $L$-function of $Y$ and~$\tY$}\label{SL}

In this section we will show that the $L$-functions
of $\Het^3(Y)$ and $\Het^3(\tY)$ are modular.
In fact, the $L$-functions are equal and belong
to the modular form
\begin{eqnarray*}
f(q)&:=&(\eta(q)\eta(q^2)\eta(q^3)\eta(q^6))^2\\
&=& q\prod_{n=1}^\infty (1-q^n)^2(1-q^{2n})^2(1-q^{3n})^2(1-q^{6n})^2\\
&=&q
-2q^2
-3q^3
+4q^4
+6q^5
+6 q^6
-16q^7\\
&&\qquad
-8q^8
+9 q^9
-12q^{10}
+12q^{11}
-12q^{12}
+38q^{13}
+\cdots\\
&=:& \sum_{n=1}^\infty a_nq^n.
\end{eqnarray*}
$f$~is a cusp form of weight~4 with respect to the groups
$$
\Gamma_0(6):=\left\{\begin{pmatrix}a & b\\ c &
d\end{pmatrix}\in\SL_2(\ZZ)\colon
c\equiv 0\bmod 6
\right\}
$$
and
$$
\Gamma_1(6):=\left\{\begin{pmatrix}a & b\\ c &
d\end{pmatrix}\in\Gamma_0(6)\colon
a\equiv 1\bmod 6
\right\}.
$$
These groups define the same modular forms
since their images in~$\PSL_2(\ZZ)$ coincide.
For the same reason, the natural map $X_1(6)\to X_0(6)$
between the associated modular curves is an isomorphism.
We denote the universal elliptic curve over~$X_1(6)$ by $S_1(6)$.
Standard formulae imply that the space $S_4(\Gamma_0(6))$
of cusp forms of weight~4 has dimension~1
and it is well known that $f$~generates this space
\cite{Shi}. Furthermore, it is a
(normalized) newform,
since there are no cusp forms of weight~4 and level~2 or~3.

The $L$-function of~$f$ is the Mellin transform
$$
L(f,s):=\sum_{n=1}^\infty a_nn^{-s}
$$
of~$f$. Since $a_n=O(n^{3/2})$ for $n\to\infty$, it converges
for $\Re(s)>5/2$. It has an analytic continuation
to an entire function. Furthermore,
there is a functional equation relating $L(f,s)$
and $L(f,4-s)$.  Since $f$~is a Hecke eigenform,
its $L$-function is an Euler product
$$
L(f,s)=\prod_{\mbox{\scriptsize $p$ prime}} L_p(f,s)
$$
with Euler factors
$$
L_p(f,s)=\frac1{1-a_pp^{-s}+p^3\cdot p^{-2s}}\qquad\mbox{for $p\ge5$}
$$
and
$L_p(f,s)=(1+p\cdot p^{-s})^{-1}$ for $p<5$.

The $L$-function of the Galois module~$\Het^3(Y)$ is also
an Euler product. Its Euler factor at a prime $p\ge5$ is
$$
L_p(\Het^3(Y),s)=\frac1{1-a_p(Y)p^{-s}+p^3\cdot p^{-2s}},
$$
where $a_p(Y):=t_3:=\Tr(\Frob_p^*;\Het^3(Y))$.
The $L$-function of~$\Het^3(\tY)$ is defined
analogously.

Our proof of the modularity of~$L(\Het^3(Y),s)$
will follow the lines of \cite[\S3.5]{V1}.
The main point is that a theorem due to Serre
(see also \cite{Sch}),
based on Faltings work and recast by Livn\'e,  essentially
reduces the proof
to checking equality of the Fourier coefficients
for only finitely many primes!
To make this precise we
recall Livn\'e's theorem   \cite[Theorem~4.3]{L}.

\begin{thm}\label{Liv} Let $S$ be a finite set of prime numbers
and let $\rho_1,\rho_2$ be  continuous
2-dimensional 2-adic representations of~$\Gal(\bQQ/\QQ)$
unramified outside~$S$.
Let $\QQ_S$ be the compositum of all quadratic extensions of~$\QQ$
which are unramified outside~$S$. Let $T$ be a set
of primes, disjoint from~$S$,
such that $\Gal(\QQ_S/\QQ)=\{\Frob_p|_{\QQ_S}\colon p\in T\}$.
Suppose that
\begin{enumerate}
\item[(a)] $\Tr\rho_1(\Frob_p)=\Tr\rho_2(\Frob_p)$ for all $p\in T$;
\item[(b)] $\det\rho_1(\Frob_p)=\det\rho_2(\Frob_p)$ for all $p\in T$;
\item[(c)] $\Tr\rho_1\equiv\Tr\rho_2\equiv 0\bmod 2$
and $\det\rho_1\equiv\det\rho_2\bmod 2$.
\end{enumerate}
Then
$\rho_1$ and~$\rho_2$ have isomorphic semisimplifications, whence
 $L(\rho_1,s)=L(\rho_2,s)$.
In particular, the good Euler factors of~$\rho_1$ and~$\rho_2$
coincide.
\end{thm}

\begin{thm}\label{th32} We have $L(\Het^3(Y),s)\circeq
L(\Het^3(\tY),s)\circeq L(f,s)$,
i.e.\ the Euler factors for ${p\ge 5}$ coincide.
\end{thm}

\begin{pf}{Proof:}{$\boxtimes$} Note that
$L(\Het^3(\tY),s)\circeq L(\Het^3(Y),s)$,
since $\tY\to Y$ is a finite map inducing
an isomorphism on~$H^{3,0}$.
So it suffices to show that
$L(\Het^3(Y),s)\circeq L(f,s)$.

We work with 2-adic cohomology (because of Livn\'{e}'s theorem).
Recall that $L(f,s)=L(\rho_f,s)$, where
 $\rho_f\colon \Gal(\bQQ/\QQ)\to \Aut(F)$
is the 2-dimensional 2-adic representation associated
to~$f$ by Deligne \cite{De}. The $L$-series $L(\Het^3(Y),s)$
also depends on the Galois module structure of~$\Het^3(Y)$,
so $L(\rho_Y,s)$ would be a more appropriate notation,
where $\rho_Y$~denotes the action of~$\Gal(\bQQ/\QQ)$ on~$\Het^3(Y)$.
We apply Livn\'{e}'s theorem to $\rho_1=\rho_Y$, $\rho_2=\rho_f$ and
$S=\{2,3\}$.
We can take  $T=\{5,7,11,13,17,19,23,73\}$,
since $\QQ_S=\QQ(\sqrt{-1},\sqrt{2},\sqrt{3})=\QQ(e^{2\pi i/24})$
and the image of $\Frob_p|_{\QQ_S}$ under
$\Gal(\QQ_S/\QQ)\stackrel{\sim}{\to} (\ZZ/24\ZZ)^*$
is simply $p$~mod~24.

We claim that the conditions (a), (b) and~(c) of
Livn\'{e}'s theorem
follow from
\begin{enumerate}
\item[$(i)$] the $p$-th Fourier coefficient of~$f$ coincides
with the trace $t_3(p)$ of~$\Frob_p$
acting on~$\Het^3(Y))$
for $p\in \{5,7,11,13,17,19,23,73\}=:T$;
\item[$(ii)$] $\det(\Frob_p; \Het^3(Y))=+p^3$ for all $p\in T$;
\item[$(iii)$] $t_3(p)$ is even for all $p\ge5$.
\end{enumerate}
Indeed, Livn\'{e}'s conditions
$(a)$ and~$(b)$ specialize to  $(i)$ and~$(ii)$ respectively,
since $\det(\rho_f(\Frob_p))=p^3$ by Deligne or by inspection
of the $p$-th Euler factor of~$L(f,s)$.
By Chebotarev's density theorem it suffices to check
condition~$(c)$ in~$\Frob_p$ for almost all~$p$.
Since for $p\ge 5$ each determinant equals $\pm p^3$, we have
$\det\rho_Y(\Frob_p)\equiv 1\equiv \det\rho_f(\Frob_p)\bmod 2$.
Finally, the evenness of~$\Tr\rho_f(\Frob_p)$ for $p\ge5$
was proved \cite[Lemma~3.12]{V1}, so condition~$(c)$
indeed reduces to condition~$(iii)$.

It remains to prove the conditions $(i)$, $(ii)$ and~$(iii)$.
Recall that
\begin{eqnarray*}
t_3 &=& (1+t_2+t_4+p^3)-\#Y(\FF_p)\\
&=& (1+50p+50p^2+p^3)-\#Y(\FF_p)
\end{eqnarray*}
(cf.~Theorem~\ref{th21} and Proposition~\ref{p24}).
By the formula $\#Y(\FF_p)=\#U(\FF_p)+50p^2+50p+20$
this boils down to
$$
t_3=p^3-19-\#U(\FF_p).
$$
Using a computer to determine $\#U(\FF_p)$
we get
the following table.

\bigskip

\begin{center}
\begin{tabular}{r|r|r|r}
$p$ & $\#U(\FF_p)$ & $\#Y(\FF_p)$ & $t_3$\\ \hline
5 & 100 & 1620 & 6\\
7 & 340 & 3160 & -16\\
11 & 1300 & 7920 & 12\\
13 & 2140 & 11260 & 38\\
17 & 5020 & 20340 & -126\\
19 & 6820 & 25840 & 20\\
23 & 11980 & 39600 & 168\\
73 & 388780 & 658900 & 218
\end{tabular}
\end{center}

\bigskip

\noindent
We wrote a straightforward Maple program
to compute these numbers.
The computation of~$\#U(\FF_{23})$ took
less than a minute on a Macintosh~G4 with a 350~MHz processor;
the computation of $\#U(\FF_{73})$ took about 80~minutes.
Using more advanced techniques one can speed up
these calculations considerably.
Note that $t_3$ indeed coincides
with the corresponding Fourier coefficient of~$f$
for these primes. This proves condition~$(i)$.

Condition~$(ii)$ follows from the observation
that $t_3\neq0$ for the 8~primes above.
Indeed, if $\alpha$ and~$\beta$ denote the eigenvalues
of Frobenius on~$\Het^3$, then $\det(\Frob_p)=\alpha\beta$
is an integer of absolute value~$p^3$
and we only have to exclude the case $\{\alpha,\beta\}=\{-p^{3/2},p^{3/2}\}$.
However, in that case we would have $t_3=\alpha+\beta=0$.

Finally, we have to check that $\#U(\FF_p)$ is even for every prime~$p\ge5$.
Consider the involution $(x_0:\cdots:x_5)\mapsto (x_0^{-1}:\cdots :x_5^{-1})$
on~$U=\{x\in N\colon x_0\cdots x_5\neq0\}$.
Its fixpoints are the 10~Segre points.
Since the number of fixpoints is even, so is $\#U(\FF_p)$.
\end{pf}

\begin{rmk}\upshape This implies
that
$a_p=(1+60p+60p^2+p^3)-\#\tY(\FF_p)$
for $p\equiv 1\bmod 4$, since $t_2:=\Tr(\Frob_p;\Het^2(\tY))=60p$
and $t_4:=\Tr(\Frob_p;\Het^4(\tY))=60p^2$
in this case.
Conversely, one can use Theorem~\ref{th32} to determine
$t_2$ and~$t_4$ on~$\tY$ for $p\equiv 3\bmod 4$.
Indeed, our proof of condition~$(**)$ also shows
that the eigenvalues of~$\Frob_p$ on~$\Het^2(\tY)$
are $\pm p$ for~$p\equiv 3\bmod 4$.
Hence $(t_2,t_4)=k(p,p^2)$ for some integer~$k$.
The inequality
$$
|(1+kp+kp^2+p^3)-\#\tY(\FF_p)|<2p^{3/2}
$$
and a computer computation of~$\#\tY(\FF_p)$
then shows that $k=40$ for $p=7,11,19,\ldots,59$.
(In particular, for these primes the variety~$\tY$ does
not satisfy condition~$(*)_p$.)
Of course, we expect that
$k=40$ for all primes~$p$ satisfying $p\equiv 3\bmod 4$, $p\neq 3$, but we
did not
try to prove this.
\end{rmk}

\section{Correspondences with relatives}

We proved in section~\ref{SL} that
$$
L(\Het^3(Y),s)\circeq L(\Het^3(\tY,s)\circeq L(f,s)
$$
where $f(q)=(\eta(q)\eta(q^2)\eta(q^3)\eta(q^6))^2$
is the normalized generator of $S_4(\Gamma_0(6))=S_4(\Gamma_1(6))$.
The Tate conjecture
then says that there should be a correspondence between $Y$
and $W=S_1(6)\times_{X_1(6)}S_1(6)$.
We shall in fact see that there is a birational equivalence
between $Y$ and~$W$ which is defined over~$\QQ$.

Recall from \cite{Be} that the universal elliptic curve
$S_1(6)\to X_1(6)\cong\PP_1$
is given by the pencil
$$
(X+Y)(Y+Z)(Z+X)+tXYZ=0.
$$

\begin{thm}\label{th41}  There exists a birational equivalence between
$S_1(6)\times_{X_1(6)}S_1(6)$ and~$Y$ which is defined over~$\QQ$.
\end{thm}

\begin{pf}{Proof:}{$\boxtimes$}
One finds immediately from the defining equations of the
Barth-Nieto quintic that
$$
x_0+x_1+x_2=-(x_3+x_4+x_5),\qquad
\frac1{x_0}+\frac1{x_1}+\frac1{x_2}=
-\left(\frac1{x_3}+\frac1{x_4}+\frac1{x_5}\right).
$$
By multiplying these two equations we obtain two pencils
of cubics, namely
$$
(x_0+x_1+x_2)\left(\frac1{x_0}+\frac1{x_1}+\frac1{x_2}\right)=
(x_3+x_4+x_5)\left(\frac1{x_3}+\frac1{x_4}+\frac1{x_5}\right)=:t'.
$$
The first (and similarly the second) pencil can also be written
in the form
$$
(x_0+x_1+x_2)(x_1x_2+x_2x_0+x_0x_1)=t'x_0x_1x_2
$$
or equivalently
$$
(x_0+x_1)(x_1+x_2)(x_2+x_0)+(1-t')x_0x_1x_2=0.
$$
Setting $t=1-t'$ this is precisely Beauville's pencil.
Hence the rational map given by
$$
(x_0:\cdots :x_5)\mapsto ((x_0:x_1:x_2),t,(x_3:x_4:x_5))
$$
defines a rational equivalence between $Y\subset\PP_5$
and $W\subset \PP_2\times\PP_1\times \PP_2$.
\end{pf}

\begin{rmk}\upshape The analogous procedure
using the equations $x_0+x_1=-(x_2+x_3+x_4+x_5)$
and $\frac1{x_0}+\frac1{x_1}
=-(\frac1{x_2}+\frac1{x_3}+\frac1{x_4}+\frac1{x_5})$
gives rise to a pencil of K3-surfaces on~$Y$.
\end{rmk}

H.~Verrill has studied in \cite{V1} a Calabi-Yau variety~$\cZ_{A_3}$
which is the smooth model of the variety~$V$ given in inhomogeneous
coordinates by
$$
V\colon (1+x+xy+xyz)(1+z+yz+xyz)=\tfrac{(t+1)^2}txyz.
$$
She has shown that $L(\Het^3(\cZ_{A_3}),s)\circeq L(f,s)$
and hence one also expects a correspondence between
$Y$ and~$\cZ_{A_3}$, resp.\ $N$ and~$V$.

\begin{thm}\label{th43} There exists a birational equivalence between $N$
and~$V$
which is defined over~$\QQ$.
\end{thm}

\begin{pf}{Proof:}{$\boxtimes$} The variety~$V$ is fibred
by a pencil of K3-surfaces.
In order to find a suitable fibration of~$Y$ we consider the pencil
$H_t=\{x_0=tx_1\}$ containing the plane $F_{01}=\{x_0=x_1=0\}$
and a residual quartic surface~$X_t$, i.e.\
$$
H_t\cap N=F_{01}\cup X_t.
$$
Combining the equations $x_0=tx_1$ and $x_0+\cdots+x_5=0$
we obtain
$$
x_1=-\tfrac1{t+1}(x_2+x_3+x_4+x_5)\qquad\mbox{and}\qquad
x_0=-\tfrac t{t+1}(x_2+x_3+x_4+x_5).
$$
Substituting this into the second defining equation
of~$N$ we get
$$
\sum_{i=0}^5x_0\cdots \hat{x}_i\cdots x_5 =
(x_2+x_3+x_4+x_5)
(\tfrac t{(t+1)^2}e_1e_3-e_4)
$$
where $e_j=e_j(x_2,x_3,x_4,x_5)$ is the $j$-th elementary
symmetric function. Hence the residual quartic~$X_t$ is given by
$$
X_t\colon e_1e_3-\tfrac{(t+1)^2}te_4=0.
$$
We consider the Cremona transformation of~$\PP_3$
given by
$$
(x_2:x_3:x_4:x_5)=(1:x:xy:xyz).
$$
Note that under this transformation
\begin{eqnarray*}
e_1 &=& 1+x+xy+xyz\\
e_3 &=& x^2y\cdot (1+z+yz+xyz)\\
e_4 &=& x^2y\cdot xyz.
\end{eqnarray*}
Comparing this with the equations defining~$V$
shows that the rational map
$(x_2:x_3:x_4:x_5)\mapsto (1:x:xy:xyz)\colon
N \dashrightarrow V$ gives the desired birational equivalence.
\end{pf}

\begin{rmk}\upshape
\begin{enumerate}
\item Combining the maps from Theorems~\ref{th41} and~\ref{th43}
one obtains a birational equivalence between $V$ and~$W$.
Saito and Yui have also found an explicit birational equivalence between
these varieties
(cf.\ \cite[p.~542]{Y}). Their map, however, differs from ours.
\item By a result of Batyrev \cite{Ba}
the two smooth Calabi-Yau varieties $Y$ and $\cZ_{A_3}$
must have the same Betti numbers and hence also the same
Hodge numbers. This agrees with our computations.
Note however, that we still need the calculation of section~\ref{S2}
to compute the Hodge numbers of~$Z$.
\end{enumerate}
\end{rmk}

The birational equivalence between $N$ and~$W$ means
that one can associate to each general pair of points $(P,Q)$
on the same fibre~$E$ of $S_1(6)\to X_1(6)$
an $H_{2,2}$-invariant Kummer surface,
resp.\ two abelian surfaces with a (1,3)-polarization and a level-2
structure. This can be made explicit in the following sense.
Let $(P,Q)\in W$. Then we can choose
homogeneous coordinates $P=(x_0:x_1:x_2)$
and $Q=(x_3:x_4:x_5)$ such that $x_0+\cdots+x_5=0$.
Then $x_0,\ldots,x_5$ define a quartic surface
$$
X=\left\{\sum_{i=0}^5 x_it_i=0\right\}\subset\PP_3
$$
where the $t_i$ are the $H_{2,2}$-invariant quartic polynomials
given in \cite[p.~189]{BN}. The surface~$X$ contains 2~sets of
16~disjoint lines and taking the branched cover along these sets of 16~lines
we obtain after blowing down the $(-1)$-curves two abelian surfaces
$A$ and~$A'$. The line bundle~$\cO_X(1)$ defines a (2,6)-polarization,
and hence also a (1,3)-polarization on $A$ and~$A'$ and the action
of the Heisenberg group~$H_{2,2}$ gives a level-2 structure.

Conversely, given a period matrix $\tau=\left(\begin{smallmatrix}
\tau_1 & \tau_2\\ \tau_2 & \tau_3\end{smallmatrix}\right)$
in Siegel upper half space we can associate to it an abelian surface~$A$
together with a (1,3)-polarization and a level-2 structure.
These data define 4~explicitely known theta functions
$\theta_0,\theta_1,\theta_2,\theta_3$ such that the map defined
by these functions embed the Kummer surface~$X$
of~$A$ as an $H_{2,2}$-invariant quartic surface
(see \cite{HNS}). Using the derivatives of
$\theta_0,\theta_1,\theta_2,\theta_3$ and Heisenberg symmetry one can,
at least in principle, compute the 16~lines
on~$X$ in~$\PP_3$. There is then a unique $H_{2,2}$-invariant
quartic of the form $\{\sum x_it_i=0\}$ containing these lines and one obtains
$X=\{\sum x_it_i=0\}$. The pair $(P,Q)\in W$ with $P=(x_0:x_1:x_2)$
and $Q=(x_3:x_4:x_5)$ is then the element of~$W$
associated to the point~$[\tau]\in \cA_{1,3}(2)$.
The latter variety has $Z$~and~$\tY$ as smooth models.

\begin{rmk}\upshape
\noindent
Threefolds which have a two dimensional Galois representation in $\Het^3$
which is known or conjectured to be associated to an elliptic modular form
of weight 4 are not easy to find.
Here we list, in order of the level~$N$ of the corresponding modular form, the
ones we are aware of.
Note in that in several cases correspondences between relatives are not
(yet) known.

\begin{enumerate}
\item $N=6$. There is a unique newform, see section \ref{SL}.
This family of threefolds contains the ones discussed in this paper, the
Barth-Nieto quintic $N$,
its double cover $\tY$ which is birational to the moduli space $\cA_{1,3}(2)$,
Verrill's threefold $\cZ_{A_3}$ and the fibred square $W$ of the
universal universal curve over $\Gamma_1(6)$.
These relatives are in correspondence with each other. Two other relatives
were found in \cite{vS}:
they are given by $\sigma_1=\alpha \sigma_5+\beta \sigma_2\sigma_3=0$ for
$(\alpha:\beta)=(1:1),\;(-2:1)$, where the $\sigma_i$ are the elementary
symmetric functions in $6$ variables.
No correspondences between any of these two and the other relatives are known.
Note that the Barth-Nieto quintic corresponds to $(1:0)$.

\item
$N=8$. There is a unique newform in $S_4(\Gamma_0(8))$.
In this family we have the desingularization $X'$ of a Siegel modular
threefold, denoted by $X$ in
\cite[p.~56]{vGN}. The variety $X$ is a complete intersection of four
quadrics in $\PP^7$
$$
\begin{array}{rcl}
Y_0^2=X_0^2+X_1^2+X_2^2+X_3^2\\
Y_1^2=X_0^2+X_1^2-X_2^2-X_3^2\\
Y_2^2=X_0^2-X_1^2+X_2^2-X_3^2\\
Y_3^2=X_0^2-X_1^2-X_2^2+X_3^2
\end{array}.
$$
(There is a misprint in the equations in \cite{vGN};
these equations agree
with those given in
the appendix of \cite{WvG}, where $X$ is denoted by $W$
and the variables $Y_1$
and $Y_2$ are interchanged.)
A relative of $X$ is the fibred square $W$ of the universal
universal curve over $\Gamma_0(8)$.
Using the moduli interpretation of $X$ and $W$, a correspondence between
them was found in \cite{EvG}.
Using the explicit equation for (an open part of) $W$
$$
U_0+U_0^{-1}+U_1+U_1^{-1}+U_2+U_2^{-1}+U_3+U_3^{-1}=0
$$
(cf.\ \cite[p.~60]{vGN}), J. Stienstra found a dominant rational map
$X\rightarrow W$ given by
$$
U_i:=(Y_i+\sqrt{2}X_i)/(Y_i-\sqrt{2}X_i)
$$
(there are unfortunately also misprints in the formulae given in
\cite[p.~60]{vGN} for this map).

\item $N=9$.  There is a unique newform in $S_4(\Gamma_0(9))$ which is
defined by a Hecke character
of the field $\QQ(\sqrt{-3})$. A threefold (a desingularization of the
intersection of two cubics
in $\PP^5$) with this L-series is given in the appendix of \cite{WvG}.
Relatives are the fibred square of the universal universal curve over
$\Gamma_0(9)$ and a product
of 3 elliptic curves,
defined over $\QQ$, with j-invariant 0. No correspondences between these
relatives seem to be known.

\item $N=12$. There is a unique newform in $S_4(\Gamma_0(12))$.
In the appendix of \cite{WvG}  a complete intersection
of a quadric and a quartic in $\PP^5$ having this L-series is given.

\item $N=21$. The threefold defined by an equation with $(\alpha:\beta)=(-3:1)$
as in the $N=6$ case has an L-series which corresponds to a newform
in $S_4(\Gamma_0(21))$, see \cite{vS}.

\item $N=25$. There is a newform in $S_4(\Gamma_0(25))$ characterized
by the first 5~Fourier coefficients
$(a_1,\ldots,a_5)$ $=$ $(1,\,1,\,7,\,-7,\,0)$. The Schoen quintic, the
3-fold in
$\PP^5$ defined by
$$
X_0^5+X_1^5+X_2^5+X_3^5+X_4^5-5X_0X_1X_2X_3X_4X_5=0
$$
has this L-series (see \cite[Proposition~5.3]{Sch}). A relative is
indicated in Remark 5.6 of that paper
(again a fibred square of a universal elliptic curve),
but no correspondence is known.

\item $N=50$. The Hirzebruch quintic 3-fold in $\PP^4$, which has 126 nodes,
has an $L$-function which corresponds to the newform in $S_4(\Gamma_0(50))$
characterized by $(a_1,\ldots,a_5)$ $=$ $(1,\,2,\,-2,\,4,\,0)$ (see
\cite{WvG}).
No relatives are (explicitly) known, but there is of course the usual
fibred square of the universal curve.

\end{enumerate}
\end{rmk}
\bigskip

\begin{rmk}\upshape
\noindent Since $Z$ is birationally isomorphic to the Siegel modular
threefold $\cA_{1,3}(2)$,
the holomorphic three form on $Z$ corresponds to a Siegel modular cusp form
of weight $3$
for the paramodular group $\Gamma_{1,3}(2)\subset \SP(4,\QQ)$.
In \cite{GH1}, this cusp form was identified as $\Delta_1^3$.
One might expect that $\Delta_1^3$ is the Saito-Kurokawa
lift of the elliptic modular form~$f$, but we do not know whether this is
indeed the case. In \cite{vGN} (see also the previous remark for $N=8$) a
similar situation occurred.
\end{rmk}

\begin{flushleft}
Klaus Hulek, Jeroen Spandaw\\
Institut f\"ur Mathematik\\
Universit\"at Hannover\\
Postfach 6009\\
D 30060 Hannover\\
Germany\\
{\tt hulek@math.uni-hannover.de\\
spandaw@math.uni-hannover.de
}

\

Bert van Geemen\\
Dipartimento di Matematica\\
Universit\`a di Pavia\\
Via Ferrata 1\\
I-27100 Pavia \\
Italy\\
{\tt geemen@dragon.ian.pv.cnr.it}

\

Duco van Straten\\
Fachbereich Mathematik, FB17\\
Johannes Gutenberg Universit\"at\\
Staudingerweg 9\\
55099 Mainz\\
Germany\\
{\tt straten@mathematik.uni-mainz.de}

\end{flushleft}

\end{document}